\documentclass{amsart}

\setbox0=\hbox{$+$}
\newdimen\plusheight
\plusheight=\ht0
\def\+{\;\lower\plusheight\hbox{$+$}\;}

\setbox0=\hbox{$-$}
\newdimen\minusheight
\minusheight=\ht0
\def\-{\;\lower\minusheight\hbox{$-$}\;}

\setbox0=\hbox{$\cdots$}
\newdimen\cdotsheight
\cdotsheight=\plusheight
\def\cds{\lower\cdotsheight\hbox{$\cdots$}}

\newcommand{\df}{\dfrac}

\newcommand{\s}{{\sigma}}

 \renewcommand{\a}{\alpha}

\renewcommand{\l}{\lambda}

\renewcommand{\t}{\varphi}

\renewcommand{\(}{\left\(}
\renewcommand{\)}{\right\)}

\renewcommand{\i}{\infty}
\renewcommand{\b}{\beta}
\renewcommand{\pmod}[1]{\,(\textup{mod}\,#1)}

\def\receivedline{\relax}
\def\dedication#1{\receivedline\vskip4pt
{\normalsize\begin{center}#1\end{center}\vskip1sp}}

\newcommand\mbgr{\mbox{BG-rank}}
\newcommand{\beqs}{\begin{equation*}}
\newcommand{\eeqs}{\end{equation*}}
\numberwithin{equation}{section}
 \theoremstyle{plain}
\newtheorem{theorem}{Theorem}[section]
\newtheorem{lemma}[theorem]{Lemma}
\newtheorem{corollary}[theorem]{Corollary}

\begin{document}

\title[ New Identities  for 7-cores with prescribed BG-rank
] {New Identities  for 7-cores with prescribed BG-rank}
\author{Alexander Berkovich}
\author{Hamza Yesilyurt}

\address{Department of Mathematics, University of Florida, 358 Little Hall,  Gainesville, FL
  32611, USA}
\address{Bilkent University, Faculty of Science, Department of Mathematics, 06800 Bilkent/Ankara, Turkey}

\email{alexb@math.ufl.edu}

\email{hamza@fen.bilkent.edu.tr}

\begin{abstract}
Let $\pi$ be a partition. BG-rank$(\pi)$ is defined as an
alternating sum of parities of parts of $\pi$ \cite{BG1}. In
\cite{BG}, Berkovich and Garvan found theta series representations
for the $t$-core generating functions $\sum_{n\geq 0}a_{t,j}(n)q^n$,
where $a_{t,j}(n)$ denotes the number of $t$-cores of $n$ with
$\mbgr=j$. In addition, they found positive $eta$-quotient
representations for odd $t$-core generating functions with extreme
values of BG-rank. In this paper we discuss representations of this
type for all $7$-cores with prescribed BG-rank. We make an essential
use of the Ramanujan modular equations of degree $7$ \cite{III} to
prove a variety of new formulas for the $7$-core generating function
\begin{equation*}
 \prod_{j\geq
1}\frac{(1-q^{7j})^7}{(1-q^j)}.
\end{equation*}
These formulas enable us to establish a number of striking
inequalities for $a_{7,j}(n)$ with $j=-1,0,1,2$ and $a_7(n)$, such
as
\begin{equation*}
a_7(2n+2) \geq 2a_7(n),\;\;a_7(4n+6)\geq 10a_7(n).
\end{equation*}

 Here $a_7(n)$
denotes a number of unrestricted $7$-cores of $n$.
Our techniques are elementary and require creative imagination only.\\
\\
\textit{`Behind every inequality there lies an identity.'} -- Basil
Gordon
\end{abstract}

\maketitle
 \dedication{Dedicated to our nephews Sam and Yu\c{s}a  }
\section{Introduction}
A partition $\pi=(\l_1,\l_2...,\l_r)$ of $n$ is a nonincreasing
sequence of positive integers that sum to $n$. The BG-rank of $\pi$
is defined as
\begin{equation}\label{BGdef}
\text{BG-rank}(\pi):=\sum_{j=1}^r(-1)^{j+1}par(\l_j),
\end{equation}
where
\begin{equation*}
par(\l_j):=\left\{ \begin{array}{ll}
         1 &  \text{if} \qquad
\l_j \equiv 1\pmod 2\\
          0 &\text{if}\qquad \l_j \equiv 0 \pmod 2\end{array} \right.
\end{equation*}
If $t$ is a positive integer, then a partition is  a $t$-core if it
has no rim hooks of length $t$ \cite{JK}. Let $\pi_{t-\text{core}}$
denote a $t$-core partition. It is shown in \cite[eq.(1.9)]{BG} that
if $t$ is odd,  then
\begin{equation}
-\Bigr\lfloor{\df{t-1}{4}}\Bigr\rfloor \leq
\text{BG-rank}(\pi_{t-\text{core}}) \leq
\Bigr\lfloor{\df{t+1}{4}}\Bigl\rfloor.
\end{equation}
Let $a_t(n)$ be the number of $t$-core partitions of $n$. It is well
known that \cite{K}, \cite{GKS}
\begin{equation}
\sum_{n \geq 0}a_t(n)q^n=\sum_{\substack{\overrightarrow{n}
 \in \mathbb{Z}^t,\;\; \overrightarrow{n}.\overrightarrow{1_t}=0}}
 q^{\tfrac{t}{2}\|\overrightarrow{n}\|^2+\overrightarrow{b_t}.\overrightarrow{n}}=\df{(q^t;q^t)_\i^t}{(q;q)_\i}=\df{E^t(q^t)}{E(q)},
\end{equation}
where
\begin{equation}
\overrightarrow{b_t}:=(0,1,2,...,t-1),\quad
\overrightarrow{1_t}:=(1,1,...,1),
\end{equation}
\begin{align*}
(a;q)_n&=(a)_n:=(1-a)(1-aq)\ldots(1-aq^{n-1}),\notag\\
(a;q)_\i&:=\prod_{n=0}^\i(1-aq^n), \quad |q|<1,\notag\\
E(q)&:=(q;q)_\i.
\end{align*}
The product $\prod_{i>0}E^{\delta_i}(q^i)$ with $\delta_i$ $\in$
$\mathbb{Z}$ will be referred to as an eta-quotient.

Next, we  recall  Ramanujan's definition for a general theta
function. Let

\begin{equation}\label{generaltheta}
f(a,b) := \sum_{n=-\i}^{\i}a^{n(n+1)/2}b^{n(n-1)/2}, \qquad |ab| <
1.
\end{equation}
The function $ f(a,b) $ satisfies the well-known Jacobi triple
product identity \cite[p.~35, Entry 19]{III}
\begin{equation}\label{19III}
f(a,b) = (-a;ab)_{\i}(-b;ab)_{\i}(ab;ab)_{\i}.
\end{equation}
Two important special cases of \eqref{generaltheta} are
\begin{equation}\label{22i}
\varphi(q) := f(q,q) = \sum_{n=-\i}^{\i}q^{n^2} =
(-q;q^2)_{\i}^2(q^2;q^2)_{\i}=\df{E^5(q^2)}{E^2(q^4)E^2(q)},
\end{equation}
and
\begin{equation}\label{22ii}
\psi(q) := f(q,q^3) = \sum_{n=-\i}^{\i}q^{2n^2-n}
=(-q;q^4)_\i(-q^3;q^4)_\i(q^4;q^4)_\i=\df{E^2(q^2)}{E(q)} .
\end{equation}

The product representations in \eqref{22i}--\eqref{22ii} are special
cases of \eqref{19III}. Also, after Ramanujan, we define
\begin{equation}\label{22iv}
\chi(q) := (-q;q^2)_{\i}.
\end{equation}

Let $a_{t,j}(n)$ be the number of $t$-core partitions of $n$ with
BG-rank=$j$ and define their generating function by
\begin{equation}\label{Cdef}
C_{t,j}(q):=\sum_{n \geq 0}a_{t,j}(n)q^n.
\end{equation} In this paper, we find representations for $C_{7,0}(q)$ and
$C_{7,1}(q)$ in terms of sums of positive eta-quotients. Such
representations for $C_{7,2}(q)$ and $C_{7,-1}(q)$ are known (see
\eqref{res1}--\eqref{res2} below). Here and throughout the
manuscript  we say that a $q$-series is positive if its power series
coefficients are nonnegative. We define $P[q]$ to be the set of all
such series. Obviously, $\t(q)$, $\psi(q)$ and $E^7(q^7)/E(q)$ $\in$
$P[q]$. In fact,  Granville and Ono  showed that \cite{GO} if $t
\geq 4$, then $a_t(n)>0$ for all $n \geq 0$. Our proofs naturally
lead us to inequalities that relate the coefficients of
$C_{7,j}(q)$, $j=0,1,-1,2$, and to equalities and inequalities for
the number of $7$-cores. The main results of this paper are
organized into two theorems whose proofs are given in sections
4 and 5.

\begin{theorem}\label{ine}
For all $n\geq 0$, we have
\begin{align}
a_7(2n+2)&\geq 2a_7(n), \label{itt}\\
a_7(4n+6)&\geq 10a_7(n),\label{it}\\
a_{7,0}(n) &\geq 9a_{7,2}(n), \label{ite} \\
a_{7,1}(n)&\geq 2a_{7,-1}(n), \label{ito}\\
a_7(28n+4r) &=5a_7(14n+2r-1),\;r=1,2,6,\label{itte}\\
 a_7(28n+4r+2)&+4a_7(7n+r-1) = 5a_7(14n+2r),\;r=2,4,5.\label{itte2}
\end{align}
\end{theorem}

By equation \eqref{evenred} below, we see that \eqref{it} and
\eqref{ite} are equivalent.

\begin{theorem}\label{mth}
\begin{align}
C_{7,1}(q)&=q\df{E(q^{28})E^3(q^{14})E(q^4)}{E(q^2)}\bigr\{\sigma(q^4)+q^2\psi(q^2)\psi(q^{14})\bigl\}.\label{mthp1}\\
C_{7,0}(q)&=\omega(q^2)\Bigr\{\psi^2(q^4)\t^2(q^{14})+q^6\psi^2(q^{28})\t^2(q^2)+q^2\df{E(q^{28})E^3(q^{14})E(q^4)}{E(q^2)}\Bigl\}\notag\\
&\;\;+q^2\psi(q^4)\psi^2(q^{14})\t^3(q^{14})+2q^4\psi^3(q^2)\psi^3(q^{14})+4q^{12}\psi^2(q^{14})\psi^3(q^{28})\t(q^2),\label{mthp2}
\end{align}
where
\begin{equation}\label{Adef}
 \omega(q):=\psi(q^{4})\t(q^{14})+q^3\psi(q^{28})\t(q^2)\;\text{and}\;\sigma(q):=\t(q)\t(q^7)+4q^2\psi(q^2)\psi(q^{14}).
\end{equation}
\end{theorem}
Observe, by \eqref{19III}, that
\begin{equation}\label{3f2}
\df{E(q^{28})E^3(q^{14})E(q^4)}{E(q^2)}=f(q^2,q^{12})f(q^4,q^{10})f(q^6,q^8)\psi(q^{14}).
\end{equation}
Therefore,  each term in \eqref{mthp1} and \eqref{mthp2} is a
product of six theta functions  which are in $P[q]$. It is
instructive to compare these representations with  those given in
\eqref{msum0}--\eqref{msum1} where for example $C_{7,1}$ is
expressed as sum of 21 multi-theta functions.

Our proofs employ the theory of modular equations. The starting
point in our proofs is one of Ramanujan's modular equations of
degree seven from which we obtained the identity
\begin{equation}\label{fdf}
\df{E^7(q^7)}{E(q)}=f(q,q^{13})f(q^3,q^{11})f(q^5,q^9)\t(q^7)\s(q^2)+8q^6\df{E^7(q^{28})}{E(q^4)}.
\end{equation}
Using several results from Ramanujan's notebooks we obtained the
following new analog of \eqref{fdf},
\begin{equation}\label{fff}
\df{E^7(q^7)}{E(q)}=f(q,q^6)f(q^2,q^5)f(q^3,q^4)\psi(q^7)\omega(q)+q^2\df{E^7(q^{14})}{E(q^2)}.
\end{equation}
The identity \eqref{fff} provided a natural compliment to
\eqref{fdf} and was essential to our proofs. For proofs of
\eqref{fdf} and \eqref{fff} see \eqref{odd22} and \eqref{iter2}.
From \eqref{fdf} and \eqref{fff}, we will deduce the following
interesting manifestly positive eta-quotient representation for the
generating function of 7-cores,
\begin{align}\label{7exp}
\df{E^7(q^7)}{E(q)} =& \s(q^4)f(q,q^{13}) f(q^3,q^{11}) f(q^5,q^9)
\t(q^7)\notag\\
&+ 2q^3\df{E^3(q^{28})E^2(q^{14})E^3(q^4)}{E^2(q^2)} +
6q^6\df{E^7(q^{28})}{E(q^4)}+2q^2\df{E(q^{14})^7}{E(q^2)}.
\end{align}
Observe, by \eqref{19III}, that

\begin{equation}\label{way}
\df{E^3(q^{28})E^2(q^{14})E^3(q^4)}{E^2(q^2)}=f(q^2,q^{12})f(q^6,q^8)f(q^4,q^{10})\psi(q^2)\psi^2(q^{14}),
\end{equation}

\begin{align}\label{hdm2}
f(q,q^{13})f(q^3,q^{11})f(q^5,q^9)\t(q^{7})=\df{\psi(q)\psi(q^7)E^4(q^{14})}{E(q^4)E(q^{28})}.
\end{align}

 The proof of \eqref{7exp} is given at the end of section 5.

In \cite{BG}, it is shown that the generating functions
$C_{t,j}(q)$, $t$ odd, can be written as sums of multi-theta
functions. We record them here for the case $t=7$. Let
\begin{align*}
B=(0,1,0,1,0,1,0),\\
\widetilde{B}=(1,0,1,0,1,0,1).
\end{align*}
and for $0\leq i \leq 6$ let $\overrightarrow{e}_i$ be the standard
unit vector in $\mathbb{Z}^7$. Then

\begin{align}
&C_{7,-1}(q)=\sum_{i=0}^{6}\sum_{\substack{\overrightarrow{n}
 \in \mathbb{Z}^7,\;\; \overrightarrow{n}.\overrightarrow{1_7}=0\\\overrightarrow{n}\equiv B+{\overrightarrow{e}}_{i}\pmod
 {2\mathbb{Z}^7}}}q^{\tfrac{7}{2}\|\overrightarrow{n}\|^2+\overrightarrow{b_7}.\overrightarrow{n}},\\
&C_{7,0}(q)=\sum_{0 \leq i_0 < i_1 < i_2 \leq
6}\sum_{\substack{\overrightarrow{n}
 \in \mathbb{Z}^7,\;\; \overrightarrow{n}.\overrightarrow{1_7}=0\\\overrightarrow{n}\equiv B
 +{\overrightarrow{e}}_{i_0}+{\overrightarrow{e}}_{i_1}+{\overrightarrow{e}}_{i_2}\pmod
 {2\mathbb{Z}^7}}}q^{\tfrac{7}{2}\|\overrightarrow{n}\|^2+\overrightarrow{b_7}.\overrightarrow{n}},\label{msum0}\\
&C_{7,1}(q)=\sum_{0 \leq i_0 < i_1  \leq
6}\sum_{\substack{\overrightarrow{n}
 \in \mathbb{Z}^7,\;\; \overrightarrow{n}.\overrightarrow{1_7}=0\\\overrightarrow{n}\equiv \widetilde{B}
 +{\overrightarrow{e}}_{i_0}+{\overrightarrow{e}}_{i_1}\pmod
 {2\mathbb{Z}^7}}}q^{\tfrac{7}{2}\|\overrightarrow{n}\|^2+\overrightarrow{b_7}.\overrightarrow{n}},\label{msum1}\\
&C_{7,2}(q)=\sum_{\substack{\overrightarrow{n}
 \in \mathbb{Z}^7,\;\; \overrightarrow{n}.\overrightarrow{1_7}=0\\\overrightarrow{n}\equiv \widetilde{B}
 \pmod
 {2\mathbb{Z}^7}}}q^{\tfrac{7}{2}\|\overrightarrow{n}\|^2+\overrightarrow{b_7}.\overrightarrow{n}}.
\end{align}

Eta-quotient representations for
\begin{equation}
C_{t,(-1)^{\tfrac{t-1}{2}}\bigr\lfloor{\tfrac{t}{4}}\bigl\rfloor
}(q)\;\;
\text{and}\;\;C_{t,(-1)^{\tfrac{t+1}{2}}\bigr\lfloor{\tfrac{t+1}{4}}\bigl\rfloor
}(q)
\end{equation}
are obtained in \cite[eq.(1.10)--(1.11)]{BG}. For $t=7$, they are as
follows
\begin{align}
C_{7,-1}(q)&=q^3\df{E^3(q^{28})E^2(q^{14})E^3(q^4)}{E^2(q^2)},\label{res1}\\
C_{7,2}(q)&=q^6\df{E^7(q^{28})}{E(q^4)}.\label{res2}
\end{align}

As we shall see next, it is easy to find eta-quotient
representations for $C_{7,0}(q)$ and $C_{7,1}(q)$ but these
representations are not manifestly positive.  Observe that if $\pi$
is a partition of $n$, then, by definition \eqref{BGdef},
\begin{equation}
\text{BG-rank}(\pi) \equiv n \pmod 2.
\end{equation}
Therefore, $C_{t,j}(q)$ is either an odd or an even function of $q$
with parity determined by the parity of $j$. In particular,
$C_{7,0}(q)$ and $C_{7,2}(q)$ are even functions of $q$ and
$C_{7,1}(q)$ and $C_{7,-1}(q)$ are odd functions of $q$. Moreover,

\begin{equation}
\sum_{n \geq 0} a_t(n)q^n
=\df{E^7(q^7)}{E(q)}=C_{7,-1}(q)+C_{7,0}(q)+C_{7,1}(q)+C_{7,2}(q).
\end{equation}
Therefore, by \eqref{res2},
\begin{align}
C_{7,0}(q)&=\text{even part of}\;
\Bigr\{\df{E^7(q^7)}{E(q)}\Bigl\}-C_{7,2}(q)\notag\\&=\df{1}{2}\Bigr\{\df{E^7(q^7)}{E(q)}+\df{E^7(-q^7)}{E(-q)}\Bigl\}
-q^6\df{E^7(q^{28})}{E(q^4)}\label{evenred}
\end{align}
and by \eqref{res1},

\begin{align}
C_{7,1}(q)&=\text{odd part of}\;
\Bigr\{\df{E^7(q^7)}{E(q)}\Bigl\}-C_{7,-1}(q)\notag\\&=\df{1}{2}\Bigr\{\df{E^7(q^7)}{E(q)}-\df{E^7(-q^7)}{E(-q)}\Bigl\}
-q^3\df{E^3(q^{28})E^2(q^{14})E^3(q^4)}{E^2(q^2)}.\label{oddred}
\end{align}

The rest of this paper is organized as follows. In the next section,
we give a brief introduction to modular equations. Then, we prove
three lemmas. In Lemma \ref{sanda}, we give several identities for
$\s(q)$ and $\omega(q)$, which were defined in \eqref{Adef}. The
identity \eqref{fff} in its equivalent form is proved in Lemma
\ref{fffl} (see \eqref{even}, \eqref{3f2} and \eqref{iter2}).  These
three lemmas are then used in sections 4 and 5
where we prove Theorems \ref{ine} and \ref{mth}.

\section{Modular Equations}
In this section, we give background information on modular
equations. For $0<k<1$, the complete elliptic integral of the first
kind $K(k)$, associated with the modulus $k$, is defined by

\begin{equation*}
K(k):=\int_0^{\pi/2}\df{d\theta}{\sqrt{1-k^2\sin^2\theta}}.
\end{equation*}

The number $ k^\prime := \sqrt{1-k^2}$ is called the
\textit{complementary modulus}. Let $ K, K^\prime, L, $ and $
L^\prime $ denote complete elliptic integrals of the first kind
associated with the moduli $ k, $ $ k^\prime , $ $ \ell, $ and $
\ell^\prime , $ respectively. Suppose that
\begin{equation}\label{i2}
 n\df{K^\prime}{K} = \df{L^\prime}{L}
\end{equation}
for some positive rational integer $ n. $ A relation between $ k $
and $ \ell $ induced by \eqref{i2} is called a {\it modular equation
of degree $ n.$ } There are several definitions of a modular
equation in the literature. For example, see the books by
R.~A.~Rankin \cite[p.~76]{rankin} and B.~Schoeneberg
\cite[pp.~141--142]{sc}. Following Ramanujan, set
$$ \a = k^2 \qquad \text{and} \qquad  \b = \ell^2. $$
We often say that $ \b $ has degree $ n $ over $ \a. $ If
\begin{equation}\label{i3}
q = \exp(-\pi K^{\prime}/K),
\end{equation}
two of the most fundamental relations in the theory of elliptic
functions are given by the formulas \cite[pp.~101--102]{III},
\begin{equation}\label{i4}
 \t^2(q) =\df{2}{\pi}K(k)
 \;\;\text{and}\;\;\a=k^2=1-\df{\t^4(-q)}{\t^4(q)}.
\end{equation}
The equation \eqref{i4} and elementary theta function identities
 make it possible
to write each modular equation as a theta function identity.
Ramanujan derived an extensive ``catalogue'' of formulas
\cite[pp.~122--124]{III} giving  the ``evaluations'' of $ E(q)$,
$\t(q)$, $\psi(q)$, and $ \chi(q)$ at various powers of the
arguments in terms of
$$ z := z_1 :=\df{2}{\pi}K(k) , \quad \a, \quad
\text{and}\quad q.$$ The evaluations that will be needed in this
paper are as follows,
\begin{align}
\t(q)&=\sqrt{z},\label{list1}\\
\t(-q)&=\sqrt{z}{(1-\a)}^{1/4},\label{list2}\\
\t(-q^2)&=\sqrt{z}{(1-\a)}^{1/8},\label{list3}\\
\psi(q)&=q^{-1/8}\sqrt{\tfrac{1}{2}z}\a^{1/8},\label{list4}\\
\psi(-q)&=q^{-1/8}\sqrt{\tfrac{1}{2}z}\{\a(1-\a)\}^{1/8},\label{list5}\\
\psi(q^2)&=2^{-1}q^{-1/4}\sqrt{z}\a^{1/4},\label{list6}\\
E(-q)&=2^{-1/6}q^{-1/24}\sqrt{z}\{\a(1-\a)\}^{1/24},\label{list65}\\
E(q^2)&=2^{-1/3}q^{-1/12}\sqrt{z}\{\a(1-\a)\}^{1/12},\label{list7}\\
\chi(-q^2)&=2^{1/3}q^{1/12}\a^{-1/12}(1-\a)^{1/24}.\label{list8}
\end{align}
We should remark that in the notation of \cite{III}, $E(q)=f(-q)$.
If $ q $ is replaced by $ q^n$, then the evaluations are given in
terms of
 $$  z_n := \df{2}{\pi}K(l), \quad \b, \quad
\text{and} \quad q^n,$$ where $ \b $ has degree $ n $ over $ \a$.

Lastly, the multiplier $ m $ of degree $ n $ is defined by
\begin{equation}\label{i5}
 m =\df{\t^2(q)}{\t^2(q^n)}=\df{z}{z_n}.
\end{equation}
The proofs of the following  modular equations of degree seven can be
found in \cite[p.~314, Entry 19(i),(iii)]{III},

\begin{align}
(\a \b)^{1/8}+\bigr\{(1-\a)(1-\b)\bigl\}^{1/8}&=1,\label{71}\\
\bigr(\tfrac{1}{2}\bigr(1+(\a
\b)^{1/8}+\bigr\{(1-\a)(1-\b)\bigl\}^{1/8}\bigl)\bigl)^{1/2}&=1-\bigr\{\a\b(1-\a)(1-\b)\bigl\}^{1/8},\label{72}\\
\Bigr(\df{(1-\b)^7}{(1-\a)}\Bigl)^{1/8}-\Bigr(\df{\b^7}{\a}\Bigl)^{1/8}=m\Bigr(\tfrac{1}{2}\bigr(1+(\a
\b)^{1/8}&+\Bigr\{(1-\a)(1-\b)\Bigl\}^{1/8}\Bigl)\bigl)^{1/2},\label{73}\\
m=\df{1-4\Bigr(\df{\b^7(1-\b)^7}{\a(1-\a)}\Bigl)^{1/24}}{\bigr\{(1-\a)(1-\b)\bigl\}^{\tfrac{1}{8}}-(\a
\b)^{\tfrac{1}{8}}},\;&\df{7}{m}=-\df{1-4\Bigr(\df{\a^7(1-\a)^7}{\b(1-\b)}\Bigl)^{1/24}}{\bigr\{(1-\a)(1-\b)\bigl\}^{\tfrac{1}{8}}-(\a
\b)^{\tfrac{1}{8}}},\label{74}\\
\Bigr(\df{(1-\b)^7}{(1-\a)}\Bigl)^{1/8}+\Bigr(\df{\b^7}{\a}\Bigl)^{1/8}&+2\Bigr(\df{\b^7(1-\b)^7}{\a(1-\a)}\Bigl)^{1/24}=\df{3+m^2}{4}.\label{75}
\end{align}

\section{Three Lemmas}
\begin{lemma}\label{sanda}
If $\s(q)$ and $\omega(q)$ are defined by \eqref{Adef}, then
\begin{align}
\s(q^2)&=\t(q)\t(q^7)-2q\psi(-q)\psi(-q^7),\label{sd12}\\
\s(q)&=\s(q^2)+2q\psi(q)\psi(q^{7}),\label{s1to2}\\
\omega^2(q)&=\psi(q)\psi(q^7)\bigr(\s(q^2)-q\psi(q)\psi(q^7)\bigl),\label{xq2}\\
\s^2(q^2)&=4q\omega^2(q)+\t^2(-q)\t^2(-q^7).\label{sd67}
\end{align}
\end{lemma}
\begin{proof}
 We start with two identities from \cite[pp.~304, 315,
 eq.~(19.1)]{III},
\begin{align}
\t(-q^2)\t(-q^{14})&=\t(-q)\t(-q^7)+2q\psi(-q)\psi(-q^7), \label{jb} \\
\psi(q)\psi(q^7)&=\psi(q^{8})\t(q^{28})+q^6\psi(q^{56})\t(q^4)+q\psi(q^2)\psi(q^{14})\label{psid}.
\end{align}
We will frequently use  \eqref{psid} in the form
\begin{align}
\psi(q)\psi(q^7)&=\omega(q^2)+q\psi(q^2)\psi(q^{14}).\label{pdd}
\end{align}
Using the well-known identity, \cite[p.~40, Entry 25
9(i),(ii)]{III}
\begin{equation*}
\t(q)=\t(q^4)+2q\psi(q^8),
\end{equation*}
it is easily verified that
\begin{align}
\t(q)\t(q^7)&=\t(q^4)\t(q^{28})+4q^8\psi(q^8)\psi(q^{56})\notag\\
&\;\;+2q\{\psi(q^{8})\t(q^{28})+q^6\psi(q^{56})\t(q^4)\}\notag\\
&=\s(q^4)+2q\omega(q^2).\label{tdd}
\end{align}

Using \eqref{pdd} and \eqref{tdd} in \eqref{jb}, we find that
\begin{equation}
\t(-q^2)\t(-q^{14})=\s(q^4)-2q\omega(q^2)+2q\omega(q^2)-2q^2\psi(q^2)\psi(q^4).
\end{equation}
Replacing $-q^2$ by $q$, we conclude that
\begin{equation}\label{sd1233}
\s(q^2)=\t(q)\t(q^7)-2q\psi(-q)\psi(-q^7),
\end{equation}
which is \eqref{sd12}. Similarly, using \eqref{pdd} and \eqref{tdd}
in \eqref{sd1233}, we arrive at
\begin{align}
\s(q^2)&=\t(q)\t(q^7)-2q\psi(-q)\psi(-q^7)\notag\\
&=\s(q^4)+2q\omega(q^2)-2q\omega(q^2)+2q^2\psi(q^2)\psi(q^{14})\notag\\
&=\s(q^4)+2q^2\psi(q^2)\psi(q^{14})\label{sd13},
\end{align}
which is \eqref{s1to2} with $q$ replaced by $q^2$. Lastly, by
\eqref{pdd}, \eqref{tdd},  and by the trivial identity
$\psi^2(q)=\psi(q^2)\t(q)$, we find that
\begin{align}
4\omega^2(q^2)&=\bigr(\psi(q)\psi(q^7)+\psi(-q)\psi(-q^7)\bigl)^2\notag\\
&=\psi^2(q)\psi^2(q^7)+\psi^2(-q)\psi^2(-q^7)+2\psi(q)\psi(q^7)\psi(-q)\psi(-q^7)\notag\\
&=\psi(q^2)\psi(q^{14})\bigr(\t(q)\t(q^7)+\t(-q)\t(-q^7)\bigl)\notag\\
&\;\;+2\bigr(\omega(q^2)+q\psi(q^2)\psi(q^{14})\bigl)\bigr(\omega(q^2)-q\psi(q^2)\psi(q^{14})\bigl)\notag\\
&=2\psi(q^2)\psi(q^{14})\s(q^4)+2\omega^2(q^2)-2q^2\psi^2(q^2)\psi^2(q^{14}),
\end{align}
from which \eqref{xq2} immediately follows.

The identity \eqref{sd67}, which is not employed in this manuscript,
was first proven in \cite{CC}. Here we provide a short new proof. By
\eqref{sd12} with $q$ replaced by $-q$, we find that
\begin{align}
&\s^2(q^2)-\t^2(-q)\t^2(-q^7)\notag\\&=
\Bigr(\s(q^2)-\t(-q)\t(-q^7)\Bigl)\Bigr(\s(q^2)+\t(-q)\t(-q^7)\Bigl)\notag\\
&=2q\psi(q)\psi(q^7)\Bigr(\s(q^2)+\t(-q)\t(-q^7)\Bigl).\label{sy11}
\end{align}
Now, by \eqref{sy11}, \eqref{xq2}, and \eqref{sd12} with $q$
replaced by $-q$, we deduce that
\begin{align*}
&\s^2(q^2)-\t^2(-q)\t^2(-q^7)-4q\omega^2(q)\notag\\
&=2q\psi(q)\psi(q^7)\Bigr(\s(q^2)+\t(-q)\t(-q^7)-2\s(q^2)+2q\psi(q)\psi(q^7)\Bigl)\\
&=0,
\end{align*}
which is \eqref{sd67}.

\end{proof}

\begin{lemma}\label{fffl}
With $\omega(q)$ defined by \eqref{Adef},
\begin{align}\label{even}
f(q,q^6)f(q^2,q^5)f(q^3,q^4)=q^2\psi^3(q^7)+\psi(q)\omega(q).
\end{align}
\end{lemma}

\begin{proof}
By \eqref{19III}, we find that
\begin{align}
f(q,q^6)f(q^2,q^5)f(q^3,q^4)=\df{(-q;q)_\i}{(-q^7;q^7)_\i}E^3(q^7)=\df{\chi(-q^7)}{\chi(-q)}E^3(q^7).\label{3f}
\end{align}

In \eqref{even}, if we replace  $q$ by $q^2$, and use \eqref{pdd},
and \eqref{3f}  with $q$ replaced by $q^2$, we are led to prove

\begin{equation}
\df{\chi(-q^{14})}{\chi(-q^2)}E^3(q^{14})=q^4\psi^3(q^{14})+\psi(q^2)\bigr\{\psi(q)\psi(q^7)-q\psi(q^2)\psi(q^{14})\bigl\}.\label{even1}
\end{equation}
Transforming \eqref{even1} by means of the evaluations given by
\eqref{list8}, \eqref{list7}, \eqref{list6} and \eqref{list4}, we
find that
\begin{align*}
&\df{1}{2}q^{-5/4}\sqrt{z_7^3}\df{\a^{1/12}\b^{1/6}(1-\b)^{7/4}}{(1-\a)^{1/24}}\\&=\df{1}{8}q^{-5/4}\sqrt{z_7^3}\b^{3/4}
+\df{1}{2}q^{-1/4}\sqrt{z_1}\a^{1/4}\bigr\{\df{1}{2}q^{-1}\sqrt{z_1z_7}(\a
\b)^{1/8}-\df{1}{4}q^{-1}\sqrt{z_1z_7}(\a \b)^{1/4}\bigl\}.
\end{align*}
Simplifying and using \eqref{i5}, we arrive at
\begin{equation}\label{even2}
4\Bigr(\df{\b^7(1-\b)^7}{\a(1-\a)}\Bigl)^{1/24}=\Bigr(\df{\b^7}{\a}\Bigl)^{1/8}+m(\a
\b)^{1/8}\bigr\{2(\a \b)^{1/8}-(\a \b)^{1/4}\bigl\}.
\end{equation}
Set $t:=(\a \b)^{1/8}$. Then, by \eqref{71}, we have
\begin{equation}
\bigr\{(1-\a)(1-\b)\bigl\}^{1/8}=1-t.
\end{equation}
The equation \eqref{even2} now takes the form
\begin{equation}\label{even3}
4\Bigr\{\df{\b (1-\b)}{t(1-t)}\Bigl\}^{1/3}=\df{\b}{t}+mt(2t-t^2).
\end{equation}
It is shown in \cite[pp.~319--320, eqs.~ (19.19), (19.21)]{III} that
\begin{align}
m&=\df{t-\b}{t(1-t)(1-t+t^2)},\label{mt1}
\end{align}
and
\begin{align}
(1-2t)m&=1-4\Bigr(\df{\b(1-\b)}{t(1-t)}\Bigl)^{1/3}.\label{mt2}
\end{align}
Using \eqref{mt2} in the left-hand side of \eqref{even3} and solving
for $m$, we obtain \eqref{mt1}. Hence, the proof of \eqref{even} is
complete.
\end{proof}

We now make several observations which will be used later. By
\eqref{even} and by \eqref{3f2} with $q^2$ replaced by $q$, we find
that

\begin{equation}\label{iter1}
\df{E(q^{14})E^3(q^7)E(q^2)}{E(q)}=q^2\psi^4(q^7)+\psi(q)\psi(q^7)\omega(q).
\end{equation}
Multiplying both sides of \eqref{iter1} by
$\df{E^4(q^7)}{E(q^2)E(q^{14})}$, we conclude that
\begin{equation}
\df{E^7(q^7)}{E(q)}=q^2\df{E^7(q^{14})}{E(q^2)}+\df{E(q^{14})E^3(q^7)E(q^2)}{E(q)}\omega(q),\label{iter2}
\end{equation}
which, by \eqref{3f2}, is equivalent to  \eqref{fff}.

We should remark that if $\b$ has degree seven over $\a$, then
$\a,\;\b$ and the \textit{multiplier} $m$ can be written as rational
functions of the parameter $t=(\a \b)^{1/8}$
\cite[pp.~316--319]{III}. This parametrization is a very efficient
tool in verifying modular equations of degree seven.

\begin{lemma}\label{lwdd}
\begin{align}\label{wdd}
\df{1}{2}\Bigr\{\df{E^7(q^7)}{E(q)}+\df{E^7(-q^7)}{E(-q)}\Bigl\}=5q^2\df{E^7(q^{14})}{E(q^2)}-4q^6\df{E^7(q^{28})}{E(q^4)}+E^3(q^2)E^3(q^{14}).
\end{align}
\end{lemma}
\begin{proof}
From \eqref{74}, we find that
\begin{equation}
7\Bigr(\df{\b^7(1-\b)^7}{\a(1-\a)}\Bigl)^{1/24}+m^2\Bigr(\df{\a^7(1-\a)^7}{\b(1-\b)}\Bigl)^{1/24}=\df{m^2+7}{4}.
\end{equation}
Upon comparison with \eqref{75}, we conclude that
\begin{equation}\label{wdd1}
5\Bigr(\df{\b^7(1-\b)^7}{\a(1-\a)}\Bigl)^{1/24}+m^2\Bigr(\df{\a^7(1-\a)^7}{\b(1-\b)}\Bigl)^{1/24}
=\Bigr(\df{(1-\b)^7}{(1-\a)}\Bigl)^{1/8}+\Bigr(\df{\b^7}{\a}\Bigl)^{1/8}+1.
\end{equation}
Transforming \eqref{wdd1} by means of the evaluations given by
\eqref{list65}, \eqref{list3} and \eqref{list4}, we find that
\begin{align}\label{wdd5}
10q^2\df{\sqrt{z}}{\sqrt{z_7}^7}\df{E^7(-q^7)}{E(-q)}+2m^2\df{\sqrt{z_7}}{\sqrt{z}^7}\df{E^7(-q)}{E(-q^7)}
=\df{\sqrt{z}}{\sqrt{z_7}^7}\df{\t^7(-q^{14})}{\t(-q^2)}+8q^6\df{\sqrt{z}}{\sqrt{z_7}^7}\df{\psi^7(q^{7})}{\psi(q)}+1.
\end{align}
Multiplying both sides of \eqref{wdd5}  by
$\df{\sqrt{z}}{\sqrt{z_7}^7}\df{(-q,q^2)_\i}{(-q^7,q^{14})_\i^7}$
and using \eqref{i5}, we obtain \eqref{wdd}.
\end{proof}

An interesting corollary of \eqref{wdd} will be given at the end of
the next section. We should add that \eqref{wdd} can be rewritten as
\begin{equation}
T_2\bigr(q^2\df{E^7(q^7)}{E(q)}\bigl) =  5q^2\df{E^7(q^7)}{E(q)}
+qE^3(q)E^3(q^7),
\end{equation}
where, the Hecke operator   $T_2$  is defined by
\begin{equation}
 T_2 \bigr( \sum  a(n) q^n \bigl)  = \sum (a(2n)+4a(n/2))q^n,
\end{equation}
with $a(n/2)$ = 0 if   n is odd.

\section{Proof of Theorem \ref{ine}} \label{terr1}
By \eqref{72} and \eqref{73}, we have
\begin{equation}\label{odd1}
\Bigr(\df{(1-\b)^7}{(1-\a)}\Bigl)^{1/8}-\Bigr(\df{\b^7}{\a}\Bigl)^{1/8}=m(1-\bigr\{\a\b(1-\a)(1-\b)\bigl\}^{1/8}).
\end{equation}
Transforming \eqref{odd1} by means of the evaluations given by
\eqref{list3}--\eqref{list5}, we find that
\begin{equation}
\df{\sqrt{z}}{\sqrt{z_7}^7}\df{\t^7(-q^{14})}{\t(-q^2)}-8\df{\sqrt{z}}{\sqrt{z_7}^7}q^6\df{\psi^7(q^{7})}{\psi(q)}
=\df{z}{z_7}\Bigr(1-\df{2}{\sqrt{z}\sqrt{z_7}}q\psi(-q)\psi(-q^7)\Bigl).
\end{equation}
Simplifying, and using \eqref{list1} and \eqref{sd12}, we conclude
that
\begin{equation}\label{odd2}
\df{\t^7(-q^{14})}{\t(-q^2)}-8q^{6}\df{\psi^7(q^{7})}{\psi(q)}
=\t^4(q^7)\Bigr\{\t(q)\t(q^7)-2q\psi(-q)\psi(-q^7)\Bigl\}=\t^4(q^7)\s(q^2).
\end{equation}
Multiplying both sides of \eqref{odd2} by
$\df{(-q,q^2)_\i}{(-q^7,q^{14})_\i^7}$, we find that
\begin{equation}\label{odd22}
\df{E^7(q^{7})}{E(q)}-8q^{6}\df{E^7(q^{28})}{E(q^4)}
=\df{\psi(q)\psi(q^7)E^4(q^{14})}{E(q^4)E(q^{28})}\s(q^2),
\end{equation}
which, by \eqref{hdm2}, is equivalent to \eqref{fdf}.

Next, by \eqref{pdd}, \eqref{s1to2}, and by \eqref{iter2}, we see
that
\begin{align}
&\text{even part of}\Bigr\{\df{E^7(q^{7})}{E(q)}\Bigl\}\notag\\
&=\df{E^4(q^{14})}{E(q^4)E(q^{28})}\omega(q^2)\s(q^2)+8q^6\df{E^7(q^{28})}{E(q^4)}\notag\\
&=\df{E^4(q^{14})}{E(q^4)E(q^{28})}\omega(q^2)\bigr(\s(q^4)+2q^2\psi(q^2)\psi(q^{14})\bigl)+8q^6\df{E^7(q^{28})}{E(q^4)}\notag\\
&=\df{E^4(q^{14})}{E(q^4)E(q^{28})}\omega(q^2)\s(q^4)+2q^2\df{E(q^{28})E^3(q^{14})E(q^4)}{E(q^2)}\omega(q^2)+8q^6\df{E^7(q^{28})}{E(q^4)}\notag\\
&=\df{E^4(q^{14})}{E(q^4)E(q^{28})}\omega(q^2)\s(q^4)+2q^2\Bigr\{\df{E^7(q^{14})}{E(q^2)}-q^4\df{E^7(q^{28})}{E(q^4)}\Bigl\}+8q^6\df{E^7(q^{28})}{E(q^4)}\notag\\
&=2q^2\df{E^7(q^{14})}{E(q^2)}+6q^6\df{E^7(q^{28})}{E(q^4)}+\df{E^4(q^{14})}{E(q^4)E(q^{28})}\omega(q^2)\s(q^4).\label{hdm}
\end{align}

Recall that we defined $P[q]$ to be the set of all $q$-series with
non-negative coefficients. Now, by \eqref{pdd} and \eqref{hdm2},
\begin{equation}\label{poss}
 \df{E^4(q^{14})}{E(q^4)E(q^{28})}\omega(q^2)=\text{even part of}
 \Bigr\{\df{\psi(q)\psi(q^7)E^4(q^{14})}{E(q^4)E(q^{28})}\Bigl\} \,\in
 P[q]
 \end{equation}

Therefore, we conclude

\begin{equation}
\df{E^7(q^7)}{E(q)}-2q^2\df{E^7(q^{14})}{E(q^2)} \in P[q],
\end{equation}
which is clearly equivalent to \eqref{itt}. Alternatively, one can
directly establish that
\begin{equation}\label{isol}
\df{E^4(q^{14})}{E(q^{4})E(q^{28})}\omega(q^2) =
f(q^4,q^{24})f^3(q^{12},q^{16}) + q^6f(q^{10},q^{18})f^3(q^2,q^{26})
\in P[q].
\end{equation}
We will not use \eqref{isol}, and so we forgo its proof.

 From
\eqref{hdm}, we have
\begin{equation}\label{hdm3}
\df{E^7(q^7)}{E(q)}=2q^2\df{E^7(q^{14})}{E(q^2)}+6q^6\df{E^7(q^{28})}{E(q^4)}+s(q),
\end{equation}
where $s(q) \in P[q]$. Iterating \eqref{hdm3}, we find that
\begin{align}
\df{E^7(q^7)}{E(q)}&=2q^2\Bigr(2q^4\df{E^7(q^{28})}{E(q^4)}+6q^{12}\df{E^7(q^{56})}{E(q^8)}+s(q^2)\Bigl)+6q^6\df{E^7(q^{28})}{E(q^4)}+s(q)\notag\\
&=10q^6\df{E^7(q^{28})}{E(q^4)}+s_1(q),
\end{align}
where $s_1(q) \in P[q]$.  This last identity clearly implies
\eqref{it}. We already remarked that, by equation \eqref{evenred},
\eqref{it} and \eqref{ite} are equivalent.

To prove \eqref{ito} we return to \eqref{odd22}.  We have by
\eqref{pdd}, \eqref{res1}, \eqref{s1to2} and by \eqref{xq2}

\begin{align}
&\text{odd part of}\Bigr\{
\df{E^7(q^{7})}{E(q)}\Bigl\}-3C_{7,-1}(q)\notag\\
&=q\df{\psi(q^2)\psi(q^{14})E^4(q^{14})}{E(q^4)E(q^{28})}\s(q^2)-3q^3\df{E^3(q^{28})E^2(q^{14})E^3(q^4)}{E^2(q^2)}\notag\\
&=q\df{E(q^4)E(q^{28})E^3(q^{14})}{E(q^2)}\Bigr\{\s(q^2)-3q^2\psi(q^2)\psi(q^{14})\Bigl\}\notag\\
&=q\df{E(q^4)E(q^{28})E^3(q^{14})}{E(q^2)}\Bigr\{\s(q^4)-q^2\psi(q^2)\psi(q^{14})\Bigl\}\notag\\
&=q\omega^2(q^2)\df{E^4(q^{14})}{E(q^4)E(q^{28})}.\label{odd226}
\end{align}
By \eqref{poss}, we see that
\begin{equation}
\text{odd part of}\Bigr\{ \df{E^7(q^{7})}{E(q)}\Bigl\}-3C_{7,-1}(q)
\,\in P[q],
\end{equation}
which, by \eqref{oddred}, is clearly equivalent to \eqref{ito}.

Lastly, we prove \eqref{itte} and \eqref{itte2}. Let $b(n)$ be
defined by
\begin{equation}
\sum_{n \geq 0}b(n)q^n=E^3(q)E^3(q^7).
\end{equation}
From \eqref{wdd} with $q^2$ replaced by $q$, we find that
\begin{equation}\label{wdd8}
\sum a_7(2n)q^{n}=5q\sum a_7(n)q^{n}-4q^3\sum a_7(n)q^{2n}+\sum
b(n)q^n.
\end{equation}
Equating the even indexed terms in both sides of \eqref{wdd8}, we
arrive at
\begin{equation}
a_7(4n)-5a_7(2n-1)=b(2n).
\end{equation}
Using Jacobi's well-known identity for $E^3(q)$ \cite[Thm. 357]{hr},
namely,
\begin{equation}
E^3(q)=\sum_{k=1}^{\i}(-1)^{k-1}(2k-1)q^{k(k-1)/2},
\end{equation}
we easily conclude that $b(n)
= 0$ if $n\equiv 2,4,5 \pmod 7$.
This observation together with \eqref{wdd8} implies \eqref{itte}.
The equation \eqref{itte2} is proved similarly by equating the odd
indexed terms in both sides \eqref{wdd8}.
\begin{corollary}
\begin{equation}
 3a_7(n-1)+b(n) \geq 0 \;\text{for all}\;\; n >0.\label{aposs}
\end{equation}
\end{corollary}
\begin{proof}

By \eqref{hdm}, we can write \eqref{wdd} in its equivalent form

\begin{align}
3q\df{E^7(q^7)}{E(q)}+E^3(q)E^3(q^7)=10q^3\df{E^7(q^{14})}{E(q^2)} +
\s(q^2)\omega(q)\df{E^4(q^7)}{E(q^2)E(q^{14})}.\label{okk}
\end{align}
By \eqref{poss}, we see that the right-hand side of \eqref{okk} is
in $P[q]$, from which \eqref{aposs} is immediate.
\end{proof}

\section{Proof of Theorem \ref{mth} and \eqref{7exp}} \label{tproof}

By \eqref{oddred}, \eqref{odd22}, \eqref{pdd}, \eqref{res1} and by
\eqref{s1to2}, we have that

\begin{align}
C_{7,1}(q)&=\text{odd part of}\Bigr\{
\df{E^7(q^{7})}{E(q)}\Bigl\}-C_{7,-1}(q)\notag\\
&=q\df{\psi(q^2)\psi(q^{14})E^4(q^{14})}{E(q^4)E(q^{28})}\s(q^2)-q^3\df{E^3(q^{28})E^2(q^{14})E^3(q^4)}{E^2(q^2)}\notag\\
&=q\df{E(q^{28})E^3(q^{14})E(q^4)}{E(q^2)}\Bigr\{\s(q^2)-q^2\psi(q^2)\psi(q^{14})\Bigl\}\notag\\
&=q\df{E(q^{28})E^3(q^{14})E(q^4)}{E(q^2)}\Bigr\{\s(q^4)+q^2\psi(q^2)\psi(q^{14})\Bigl\}.\label{c71}
\end{align}

This  completes the proof of \eqref{mthp1}.

Next, we prove \eqref{mthp2}. Combining \eqref{iter1} and
\eqref{iter2}, we have
\begin{equation}\label{iter3}
\df{E^7(q^7)}{E(q)}=q^2\df{E^7(q^{14})}{E(q^2)}+q^2\psi^4(q^7)\omega(q)+\psi(q)\psi(q^7)\omega^2(q).
\end{equation}

Using \eqref{iter2} with $q$ replaced by $q^2$ in  \eqref{iter3}, we
find that
\begin{align}
\df{E^7(q^7)}{E(q)}&=q^2\Bigr\{q^4\df{E^7(q^{28})}{E(q^4)}+\df{E(q^{28})E^3(q^{14})E(q^4)}{E(q^2)}\omega(q^2)\Bigl\}
+q^2\psi^4(q^7)\omega(q)+\psi(q)\psi(q^7)\omega^2(q)\notag\\
&=q^6\df{E^7(q^{28})}{E(q^4)}+q^2\df{E(q^{28})E^3(q^{14})E(q^4)}{E(q^2)}\omega(q^2)+q^2\psi^4(q^7)\omega(q)+\psi(q)\psi(q^7)\omega^2(q).\label{iter4}
\end{align}

It now remains to find the even part of the last two terms on the
right side of \eqref{iter4}. This is is easily done with the
even-odd dissections of $\omega(q)$ and $\psi(q)\psi(q^7)$ given by
\eqref{Adef} and \eqref{pdd} and  the formula (see \cite[p.~40,
Entry 25 (iv)--(vii)]{III})

\begin{equation}
\psi^4(q)=\psi^2(q^2)(\t^2(q^2)+4q\psi^2(q^4))\label{4p4}
\end{equation}
with $q$ replaced by $q^7$.

Lastly, we prove \eqref{7exp}. Arguing as in \eqref{odd226}, we find
that
\begin{align}
&\text{odd part of}\Bigr\{
\df{E^7(q^{7})}{E(q)}\Bigl\}-2C_{7,-1}(q)\notag\\
&=q\df{E(q^{28})E^3(q^{14})E(q^4)}{E(q^2)}\Bigr\{\s(q^2)-2q^2\psi(q^2)\psi(q^{14})\Bigl\}\notag\\
&=q\df{E(q^{28})E^3(q^{14})E(q^4)}{E(q^2)}\s(q^4),\label{aodd}
\end{align}
where in the last step, we used \eqref{s1to2}. Using \eqref{aodd}
together with \eqref{hdm}, and by \eqref{pdd} and \eqref{hdm2}, we
arrive at
\begin{align}
\df{E^7(q^{7})}{E(q)}&=q\df{E(q^{28})E^3(q^{14})E(q^4)}{E(q^2)}\s(q^4)+2C_{7,-1}(q)+2q^2\df{E^7(q^{14})}{E(q^2)}\notag\\
&\;\;+6q^6\df{E^7(q^{28})}{E(q^4)}+\df{E^4(q^{14})}{E(q^4)E(q^{28})}\omega(q^2)\s(q^4)\notag\\
&=2C_{7,-1}(q)+2q^2\df{E^7(q^{14})}{E(q^2)}+6q^6\df{E^7(q^{28})}{E(q^4)}
+\df{E^4(q^{14})}{E(q^4)E(q^{28})}\s(q^4)\Bigr\{\omega(q^2)+q\psi(q^2)\psi(q^{14})\Bigl\}\notag\\
&=2C_{7,-1}(q)+2q^2\df{E^7(q^{14})}{E(q^2)}+6q^6\df{E^7(q^{28})}{E(q^4)}+\df{E^4(q^{14})}{E(q^4)E(q^{28})}\s(q^4)\psi(q)\psi(q^7)\notag\\
&=2C_{7,-1}(q)+2q^2\df{E^7(q^{14})}{E(q^2)}+6q^6\df{E^7(q^{28})}{E(q^4)}+\s(q^4)f(q,q^{13})f(q^3,q^{11})f(q^5,q^{9})\t(q^7),
\end{align}
which, by \eqref{res1},  is equal to the right hand side of
\eqref{7exp}.

\section{Concluding Remarks}
The  inequalities , \eqref{itt} and \eqref{it} (or equivalently
\eqref{ite}), of Theorem \ref{ine} are not optimal. Numerical
evidence suggest that
\begin{align*}
a_7(2n+2) \geq 3a_7(n) \;\;\text{for all}\;\; n \geq 1,\\
a_7(4n+6)\geq 15a_7(n) \;\;\text{for all}\;\; n\geq 1,\\
a_7(4n+6)\geq 11a_7(n) \;\;\text{for all}\;\; n\geq 0.
\end{align*}
 Our attempts to improve
Theorems \ref{ine} and \ref{mth} led us to the following interesting
conjectures:
\begin{align}
\psi(q)\bigr(\psi^2(q)-\psi^2(q^7)\bigl)\in P[q],\\
\psi(q)\bigr(\t^2(q)-\t^2(q^7)\bigl) \in P[q],\\
 \t(q)\bigr(\psi^2(q)-\psi^2(q^7)\bigl) \in P[q],
 \end{align}
and
\begin{align}
 \psi(q)\bigr(\t^2(q)-\psi^2(q^7)\bigl) \in P[q].
\end{align}

Referee pointed out that \eqref{itte} and \eqref{itte2} extend easily using our arguments
to a few other
arithmetic progressions; for example,
\begin{equation*}
a_7(196n+4r) = 5a_7(98n+2r-1)\;\;\text{for}\; r =10, 17, 45.
\end{equation*}
\section{Acknowledgment}
We would like to thank George Andrews, Bruce Berndt, Frank Garvan
and Li-Chien Shen for their interest and helpful comments. Frank
Garvan communicated to us elegant alternative proofs of \eqref{itt},
\eqref{itte} and \eqref{aposs}.

\end{document}